\documentclass[12pt]{article}
\usepackage{amssymb,amsmath,latexsym}

\def\proof{{\it Proof}.\ }

\def\wbull{\hfill\vrule height .9ex width .8ex depth -.1ex}

\newtheorem{formula}{}[section]
\newtheorem{proposition}[formula]{Proposition}
\newtheorem{definition}[formula]{Definition}
\newtheorem{corollary}[formula]{Corollary}

\newtheorem{lemma}[formula]{Lemma}
\newtheorem{theorem}[formula]{Theorem}

\def\thrm{\begin{theorem}}
\def\thrml#1{\begin{theorem}\label{#1}}
\def\ethrm{\end{theorem}}
\def\prp{\begin{proposition}}
\def\prp#1{\begin{proposition}\label{#1}}
\def\eprp{\end{proposition}}
\def\dfntn{\begin{definition}}
\def\dfntnl#1{\begin{definition}\label{#1}}
\def\edfntn{\end{definition}}
\def\nmrt{\begin{enumerate}}
\def\enmrt{\end{enumerate}}

\def\qtn{\begin{equation}}
\def\qtnl#1{\begin{equation}\label{#1}}
\def\eqtn{\end{equation}}
\def\lmm{\begin{lemma}}
\def\lmml#1{\begin{lemma}\label{#1}}
\def\elmm{\end{lemma}}
\def\crllr{\begin{corollary}}
\def\crllrl#1{\begin{corollary}\label{#1}}
\def\ecrllr{\end{corollary}}
\def\css{\begin{cases}}
\def\ecss{\end{cases}}

\title{ \bf{On a Koolen -- Park inequality  and Terwilliger graphs}}
\author{
Alexander Gavrilyuk
\thanks{
Partially supported by RFFI grant (project no. 08-01-00009).}\\[-1pt]
\small Ural Division of the Russian Academy of Sciences\\[-3pt]
\small Institute of Mathematics and Mechanics\\[-3pt]
\small 16, S.Kovalevskaja street, 620219, Yekaterinburg, Russia\\[-3pt]
{\tt \small alexander.gavriliouk@gmail.com}\\[-3pt]}

\begin{document}

\maketitle

\begin{abstract}
J.H. Koolen and J. Park have proved a lower bound 
for intersection number $c_2$ of a distance-regular graph $\Gamma$. 
Moreover, they showed that the graph $\Gamma$ which attains the equality 
in this bound is a Terwilliger graph. We prove that $\Gamma$ is the icosahedron, 
the Doro graph or the Conway-Smith graph, if equality is attained and $c_2\ge 2$.
\medskip

\end{abstract}
\newpage

\section{Introduction}

Let $\Gamma$ be a distance-regular graph with degree $k$ and diameter at least 2.
Let $c$ be maximal such that for each vertex $x\in \Gamma$ and every pair of nonadjacent 
vertices $y,z$ of $\Gamma_1(x)$, there exists a $c$-coclique in $\Gamma_1(x)$ 
containing $y,z$. In \cite{KP}, J.H. Koolen and J. Park have shown that the following bound holds:
\qtnl{E11}
  c_2-1\ge {\rm max}\{\displaystyle {{\frac {c' (a_1+1)-k}{{c'\choose 2}}}}\ |\ 2\le c'\le c\}, 
\eqtn
and equality implies $\Gamma$ is a Terwilliger graph. 
(For definitions see sections 2 and 3.)

The similar inequality for a distance-regular graph with $c$-claw was proved 
by C.D. Godsil, see \cite{CDG}. J.H. Koolen and J. Park \cite{KP} have noted that the bound (\ref{E11})
is met exactly for all known examples of Terwilliger graphs. We recall that only three 
examples of distance-regular Terwilliger graphs with $c_2\ge 2$ are known: 
the icosahedron, the Doro graph and the Conway-Smith graph.

In this paper, we will show that the distance-regular graph $\Gamma$ with $c_2\ge 2$ 
which attains the equality in (\ref{E11}) is a known Terwilliger graph.

\section{Definitions and preliminaries}

We consider only finite, undirected graphs without loops or multiple edges.
Let $\Gamma$ be a connected graph. The \emph{distance} ${\rm d}(u,w)$ between any 
two vertices $u$ and $w$ of $\Gamma$ is the length of a shortest path from $u$ to 
$w$ in $\Gamma$. The \emph{diameter} ${\rm diam}(\Gamma)$ of $\Gamma$ is the maximal 
distance occuring in $\Gamma$.

For a subset $A$ of the vertex set of $\Gamma$, we will also write $A$ for 
the subgraph of $\Gamma$ induced by $A$. 
For a vertex $u$ of $\Gamma$, define $\Gamma_i(u)$ to be the set 
of vertices which are at distance $i$ from $u$ ($0\le i\le {\rm diam}(\Gamma))$.
The subgraph $\Gamma_1(u)$ is called the \emph{local graph} of a vertex $u$ 
and the \emph{degree} of $u$ is the number of neighbours of $u$, i.e. $|\Gamma_1(u)|$.

For two vertices $u,w\in \Gamma$ with ${\rm d}(u,w)=2$, the subgraph 
$\Gamma_1(u)\cap \Gamma_1(w)$
is called $\mu$-\emph{subgraph} of vertices $u,w$. We say the number $\mu(\Gamma)$ is 
\emph{well-defined}, if each $\mu$-subgraph occuring in $\Gamma$ contains the same number 
of vertices which is equal to $\mu(\Gamma)$.
\medskip

Let $\Delta$ be a graph. 
A graph $\Gamma$ is \emph{locally} $\Delta$, 
if, for all $u\in \Gamma$, the subgraph $\Gamma_1(u)$ is isomorphic 
to $\Delta$. 
A graph is \emph{regular} with degree $k$, 
if the degree of each its vertex is $k$.
\medskip

A connected graph $\Gamma$ with diameter $d={\rm diam}(\Gamma)$ is \emph{distance-regular}, 
if there are integers $b_i$, $c_i$ ($0\le i\le d$) such that for any two vertices $u,w\in \Gamma$ 
with ${\rm d}(u, w)=i$, there are exactly $c_i$ neighbours of $w$ in $\Gamma_{i-1}(u)$ and $b_i$ 
neighbours of $w$ in $\Gamma_{i+1}(u)$ 
(we assume that $\Gamma_{-1}(u)$ and $\Gamma_{d+1}(u)$ are empty sets). 
In particular, distance-regular graph $\Gamma$ is
regular with degree $b_0$, $c_1=1$ and $c_2=\mu(\Gamma)$. 
For each vertex $u\in \Gamma$ and $0\le i\le d$, the subgraph $\Gamma_i(u)$ is regular with 
degree $a_i=b_0-b_i-c_i$. The numbers $b_i$, $c_i$ ($0\le i\le d$) are called the \emph{intersection 
numbers} and the array $\{b_0,b_1,\ldots,b_{d-1};c_1,c_2,\ldots,c_{d}\},$ is called 
the \emph{intersection array} of distance-regular graph $\Gamma$.

\medskip

A graph $\Gamma$ is \emph{amply regular} with parameters $(v,k,\lambda,\mu)$, 
if $\Gamma$ has $v$ vertices, it is regular with degree $k$ and the following two 
conditions hold:
 
\quad $i$) for each pair of adjacent vertices $u,w\in \Gamma$, the subgraph 
$\Gamma_1(u)\cap \Gamma_1(w)$ contains exactly $\lambda$ vertices; 

\quad $ii$) $\mu=\mu(\Gamma)$ is well-defined. 
\medskip

An amply regular graph with diameter 2 is called a \emph{strongly regular} graph and 
it is a distance-regular graph. A distance-regular graph is an amply regular graph 
with parameters $k=b_0$, $\lambda=b_0-b_1-1$ and $\mu=c_2$.
\medskip

Recall that a ($c-$)clique (or \emph{complete} graph) is a graph (on $c$ vertices) 
in which every pair of its vertices is adjacent. A ($c-$)coclique is a graph (on $c$ vertices) 
in which every pair of its vertices is not adjacent. 
\medskip

Let $\Gamma$ be a strongly regular graph with parameters $(v,k,\lambda,1)$. 
There are integers $r$ and $s$ such that the local graph of each vertex of $\Gamma$ 
is the disjoint union of $r$ copies of $s$-clique. Futhermore, $v=1+rs+s^2r(r-1)$, 
$k=rs$ and $\lambda=s-1$. We denote the set of strongly regular graphs with such 
parameters by ${\cal F}(s,r)$.

Any graph of ${\cal F}(1,r)$, i.e. a strongly regular graph with $\lambda=0$ and $\mu=1$, 
is called a \emph{Moore} strongly regular graph. It is well known (see Chapter 1 \cite{BCN})
that any Moore strongly regular graph has degree 2, 3, 7 or 57. 
The graphs with degree 2, 3 and 7 are the pentagon, the Petersen graph and 
the Hoffman Ц- Singleton graph, respectively. Whether a Moore graph with degree 57 
exists is an open problem.

\lmml{L21}
Suppose that ${\cal F}(s,r)$ is nonempty set of graphs. 
Then $s+1\le r$.
\elmm
\proof Let $\Gamma$ be a graph of ${\cal F}(s,r)$. 
We may choose vertices $u$ and $w$ of $\Gamma$ with ${\rm d}(u,w)=2$. 
Let $x$ be a vertex of $\Gamma_1(u)\cap \Gamma_1(w)$. Then the subgraph  
$\Gamma_1(w)-(\Gamma_1(x)\cup \{x\})$ contains a coclique of size at most 
$r-1$. Let us consider a $s$-clique of $\Gamma_1(u)-\Gamma_1(w)$
on vertices $y_1,y_2,..,y_s$. The subgraph $\Gamma_1(w)\cap \Gamma_1(y_i)$ 
($1\le i\le s$) contains a single vertex $z_i$. 
The vertices $z_1,z_2,..,z_s$ are mutually nonadjacent and distinct.
Hence, $s\le r-1$. The lemma is proved.\wbull

\section{Terwilliger graphs}

In this section we give a definition of Terwilliger graphs and 
some useful facts concerning them.
\medskip

A \emph{Terwilliger graph} is a connected noncomplete graph $\Gamma$ such 
that $\mu(\Gamma)$ is well-defined and each $\mu$-subgraph occuring in $\Gamma$ 
is a complete graph (hence, there are no induced quadrangles in $\Gamma$). 
If $\mu(\Gamma)>1$, then, for each vertex $u\in \Gamma$, the local graph 
of $u$ will also be a Terwilliger graph with diameter 2 and 
$\mu(\Gamma_1(u))=\mu(\Gamma)-1$.

For an integer $\alpha\ge 1$, a $\alpha$-clique extension of a graph $\bar \Gamma$ 
is the graph $\Gamma$ obtained from $\bar \Gamma$ by replacing each vertex 
$\bar u\in \bar \Gamma$ by a clique $U$ of $\alpha$ vertices, where for any 
$\bar u,\bar w\in \bar \Gamma$, $u\in U$ and $w\in W$,
$\bar u$ and $\bar w$ are adjacent if and only if $u$ and $w$ are adjacent.

\lmml{L31}
Let $\Gamma$ be an amply regular Terwilliger graph with parameters
$(v,k,\lambda,\mu)$, where $\mu>1$. There is the number $\alpha$ such that 
the local graph of each its vertex is the $\alpha$-clique extension of a 
strongly regular Terwilliger graph with parameters $(\bar v, \bar k,\bar \lambda,\bar \mu)$, 
where
$$\bar v=k/\alpha,\ \bar k=(\lambda-\alpha+1)/\alpha,\ \bar \mu=(\mu-1)/\alpha,$$
and $\alpha\le \bar \lambda+1$. In particular, if $\bar \lambda=0$, then $\alpha=1$.
\elmm
\proof The result follows from \cite[Theorem~1.16.3]{BCN}.\wbull
\medskip

We know only three examples of amply regular Terwilliger graphs with $\mu\ge 2$.
All of them are unique distance-regular locally Moore graphs:

$(1)$ the icosahedron with intersection array $\{5,2,1;1,2,5\}$ is locally pentagon.

$(2)$ the Doro graph with intersection array $\{10,6,4;1,2,5\}$ and the Conway-Smith 
graph with intersection array $\{10,6,4,1;1,2,6,10\}$ are locally Petersen graphs. 
\medskip

In \cite{GM}, A. Gavrilyuk and A. Makhnev have shown that a distance-regular 
locally Hoffman -- Singleton graph has intersection array $\{50,42,9;1,2,42\}$ 
or $\{50,42,1;1,2,50\}$ and hence it is a Terwilliger graph.
Whether the graphs with these intersection arrays exist is an open question.

\lmml{L32}
Let $\Gamma$ be a Terwilliger graph. Suppose that, 
for an integer $\alpha\ge 1$, the local graph of each its vertex is 
the $\alpha$-clique extension of a Moore strongly regular graph $\Delta$.
Then $\alpha=1$ and one of the following holds:

$(1)$ $\Delta$ is the pentagon and $\Gamma$ is the icosahedron$;$

$(2)$ $\Delta$ is the Petersen graph and $\Gamma$ is the Doro graph 
or the Conway-Smith graph$;$

$(3)$ $\Delta$ is the Hoffman -- Singleton graph or a graph with degree $57$, 
in both cases diameter of $\Gamma$ is at least $3$.
\elmm
\proof It is easy to see that the graph $\Gamma$ is amply regular. By Lemma \ref{L31}, 
we have $\alpha=1$. The statements (1) and (2) follow from \cite[Proposition~1.1.4]{BCN} 
and \cite[Theorem~1.16.5]{BCN}, respectively. 

If the graph $\Delta$ is the Hoffman -- Singleton graph and diameter of $\Gamma$ is 2, 
then $\Gamma$ is strongly regular with parameters $(v,k,\lambda,\mu)$, where $k=50$, 
$\lambda=7$ and $\mu=2$. By \cite[Theorem~1.3.1]{BCN}, the eigenvalues of $\Gamma$ are $k$ 
and the roots of the quadratic equation $x^2+(\mu-\lambda)x+(\mu-k)=0$. The roots 
of the equation $x^2-5x-48=0$ are not integers, that is impossible.
In the remained case, when $\Delta$ is regular with degree 57, we will get the same 
contradiction. The lemma is proved.\wbull
\medskip

The next lemma is useful in the proof of Theorem \ref{T41} (see Section 4).

\lmml{L33}
Let $\Gamma$ be a strongly regular Terwilliger graph with parameters 
$(v,k,\lambda,\mu)$. Suppose that, for an integer $\alpha\ge 1$, the local graph of each its 
vertex is the $\alpha$-clique extension of a strongly regular graph $\Delta$ with parameters 
$(\bar v,\bar k,\bar \lambda,\bar \mu)$. Then $\bar k-\bar \lambda-\bar \mu>1$ implies that 
$k-\lambda-\mu>1$.
\elmm
\proof We have $k=\alpha(1+\bar k+\bar k(\bar k-\bar \lambda-1)/\bar \mu)$, 
$\lambda=\alpha \bar k+\alpha-1$ and $\mu=\alpha\bar \mu+1$. If 
$\bar k-\bar \lambda-\bar \mu>1$, then $\bar k(\bar k-\bar \lambda-1)/\bar \mu>\bar k$
and this implies that 
$k-\lambda-\mu=\alpha(\bar k(\bar k-\bar \lambda-1)/\bar \mu-\bar \mu)>\alpha(\bar k-\bar \mu)>\alpha(\bar \lambda+1)\ge 1$.\wbull

\section{Koolen -- Park inequality}
In this section, we consider the bound (\ref{E11}) and classify distance-regular 
graphs with $c_2\ge 2$ which attain this bound. 
\medskip

The next proposition is a slight generalization of \cite[Proposition~3]{KP}. J.H. Koolen 
and J. Park \cite[Proposition~3]{KP} formulated the next proposition for distance-regular graphs. 
We generalize it to amply regular graphs. (Our proof is similar to the one in 
J.H. Koolen and J. Park \cite{KP}, but we give it for convenience of the reader.)

\prp{P41}
Let $\Gamma$ be an amply regular graph with parameters 
$(v,k,\lambda,\mu)$ and $c\ge 2$ be maximal such that for each vertex $x\in \Gamma$ 
and every pair of nonadjacent vertices $y,z$ of $\Gamma_1(x)$,
there exists a $c$-coclique in $\Gamma_1(x)$ containing $y,z$. Then 
$$
  \mu-1\ge {\rm max}\{\displaystyle {{\frac {c' (\lambda+1)-k}{{c'\choose 2}}}}\ |\ 2\le c'\le c\},
$$ 
and if equality is attained, then $\Gamma$ is a Terwilliger graph.
\eprp
\proof Let $\Gamma_1(x)$ contain a coclique $C'$ on vertices 
$y_1,y_2,\ldots,y_{c'}$, $c'\ge 2$. 
Since ${\rm d}(y_i,y_j)=2$, $|\Gamma_1(x)\cap \Gamma_1(y_i)\cap \Gamma_1(y_j)|\le \mu-1$ 
holds for all $i\ne j$. Then by the principle of inclusion and exclusion,
$$k=|\Gamma_1(x)|\ge |\cup_{i=1}^{c'}(\Gamma_1(x)\cap (\Gamma_1(y_i)\cup \{y_i\}))|$$
$$\ge \sum_{i=1}^{c'}|\Gamma_1(x)\cap (\Gamma_1(y_i)\cup \{y_i\})|-\sum_{1\le i<j\le c'}|\Gamma_1(x)\cap \Gamma_1(y_i)\cap \Gamma_1(y_j)|$$
$$\ge c'(\lambda+1)-{c' \choose 2}(\mu-1).$$ 

So,
\qtnl{E41}
\mu-1 \ge \displaystyle {{\frac {c' (\lambda+1)-k}{{c'\choose 2}}}}. 
\eqtn

Note that equality in (\ref{E41}) implies that 
$\Gamma_1(x)\subseteq \cup_{i=1}^{c'}(\Gamma_1(y_i)\cup \{y_i\})$ 
holds and we have $|\Gamma_1(x)\cap \Gamma_1(y_i)\cap \Gamma_1(y_j)|=\mu-1$ 
for all $i\ne j$.

Let $c$ be maximal satisfying the condition of the Proposition \ref{P41}. Then
\qtnl{E42}
\mu-1 \ge {\rm max}\{\displaystyle {{\frac {c' (\lambda+1)-k}{{c'\choose 2}}}}\ |\ 2\le c'\le c\}.
\eqtn

We may assume that for an integer $c''$, where $2\le c''\le c$,
equality holds in (\ref{E42}), i.e. 
\qtnl{E43}
\mu-1={\frac {c'' (\lambda+1)-k}{{c''\choose 2}}}={\rm max}\{\displaystyle {{\frac {c' (\lambda+1)-k}{{c'\choose 2}}}}\ |\ 2\le c'\le c\}.
\eqtn

We will show $c=c''$. 
For a vertex $x\in \Gamma$ and nonadjacent vertices $y,z\in \Gamma_1(x)$, 
there exists a $c$-coclique $C$ in $\Gamma_1(x)$ containing $y,z$. 
The equality (\ref{E43}) implies that, for any subset of 
vertices $\{y_1,y_2,\ldots,y_{c''}\}\subseteq C$, 
$\Gamma_1(x)\subseteq \cup_{i=1}^{c''}(\Gamma_1(y_i)\cup \{y_i\})$ holds.
But if $c''<c$, then $C\not\subset \cup_{i=1}^{c''}(\Gamma_1(y_i)\cup \{y_i\})$, 
which is the contradiction.

Hence, $c=c''$ and we have 
$|\Gamma_1(x)\cap \Gamma_1(y)\cap \Gamma_1(z)|=\mu-1$ for every pair of nonadjacent 
vertices $y,z\in \Gamma_1(x)$ and for all $x\in \Gamma$. 
This implies that each $\mu$-subgraph occuring in $\Gamma$ is a clique of size $\mu$ 
and $\Gamma$ is a Terwilliger graph.\wbull
\medskip

We call the inequality (\ref{E42}) $\mu$-\emph{bound}.
\medskip

Let $\Gamma$ be an amply regular Terwilliger graph with parameters $(v,k,\lambda,\mu)$. 
If $\mu=1$, then the local graph of each its vertex is the 
disjoint union of $k/(\lambda+1)$ copies of $(\lambda+1)$-clique, so 
equality in $\mu$-bound is attained. If $\mu\ge 2$, then we know only three examples 
of $\Gamma$ (see Section 3) with $\mu=2$ and each of them attains equality in $\mu$-bound:

$(1)$ $\Gamma$ is the icosahedron. The pentagon contains a 2-coclique and is regular with 
degree 2, i.e. $c=2$ and $\lambda=2$, hence we have $(2\cdot (2+1)-5)/{2 \choose 2}=1=\mu-1$.

$(2)$ $\Gamma$ is the Doro graph or the Conway-Smith graph. The Petersen graph contains 
a 4-coclique and is regular with degree 3, hence we have 
$(4\cdot (3+1)-10)/{4 \choose 2}=(16-10)/6=1=\mu-1$.
\medskip

Recall that the Hoffman -- Singleton graph contains a 15-coclique. 
If $\Gamma$ is an amply regular locally Hoffman -- Singleton graph and 
is a Terwilliger graph, then $\mu=2$, but equality in $\mu$-bound is not attained. 

\thrml{T41}
Let $\Gamma$ be an amply regular graph with parameters 
$(v,k,\lambda,\mu)$ and $\mu>1$. If $\Gamma$ attains equality in $\mu$-bound, then 
$\mu=2$ and $\Gamma$ is the icosahedron, the Doro graph or the Conway-Smith graph.
\ethrm
\proof By Proposition \ref{P41}, the graph $\Gamma$ is a Terwilliger graph and, 
be Lemma \ref{L31}, there is an integer $\alpha\ge 1$ such that the local graph 
of each its vertex is the $\alpha$-clique extension of a strongly regular 
Terwilliger graph with parameters $(\bar v,\bar k,\bar \lambda,\bar \mu)$.
By Lemma \ref{L31}, we have $k=\alpha \bar v$, $\lambda=\alpha \bar k+(\alpha-1)$ 
and $\mu=\alpha \bar \mu+1$. 

By the assumption on $\Gamma$, for a vertex $u\in \Gamma$, the local graph of $u$ 
contains a $c$-coclique which attains equality in $\mu$-bound, i.e.  
$$\mu(\Gamma)-1=\alpha \bar \mu={\frac {c (\alpha \bar k+(\alpha-1)+1)-\alpha \bar v}{{c\choose 2}}}=\alpha {\frac {c (\bar k+1)-\bar v}{{c\choose 2}}}$$
and
$$\bar \mu={\frac {c (\bar k+1)-\bar v}{{c\choose 2}}}.$$

Straightforward,
$$c^2\bar \mu-c(\bar \mu+2(\bar k+1))+2\bar v=0,$$
so 
$$c={\frac {(\bar \mu+2(\bar k+1))\pm \sqrt{(\bar \mu+2(\bar k+1))^2-8\bar v\bar \mu}}{2\bar \mu}},$$ 
and 
$$(\bar \mu+2(\bar k+1))^2\ge 8\bar v\bar \mu.$$

Let the subgraph $\Gamma_1(u)$ be isomorphic to the $\alpha$-clique extension of 
a strongly regular Terwilliger graph $\Delta$. The cardinality of the 
vertex set of $\Delta$ is equal to $\bar v=1+\bar k+\bar k(\bar k-\bar \lambda-1)/\bar \mu$, 
hence:
$$(\bar \mu+2(\bar k+1))^2\ge 8(\bar \mu+\bar k\bar \mu+\bar k(\bar k-\bar \lambda-1)),$$
$$\bar \mu^2+4\ge 4\bar \mu+4\bar k\bar \mu+4\bar k^2-8\bar k\bar \lambda-16\bar k.$$

Next, 
$$(\bar \mu/2)^2+1\ge \bar \mu+\bar k\bar \mu+\bar k^2-2\bar k\bar \lambda-4\bar k,$$
\qtnl{E44}
((\bar \mu/2)-(\bar k+1))^2\ge 2\bar k(\bar k-\bar \lambda-1).
\eqtn

At first, we may assume $\bar \mu=1$. There are integers $s,r$ such that 
$\Delta\in {\cal F}(s,r)$ and $\bar k=rs$, $\bar \lambda=s-1$.
If $\bar k-\bar \lambda-1\ge \bar k/2+1$, then 
$2\bar k(\bar k-\bar \lambda-1)\ge 2\bar k(\bar k/2+1)=\bar k^2+2\bar k$. 
It follows from (\ref{E44}) that $(\bar k+1/2)^2\ge \bar k^2+2\bar k$ and hence $1/4\ge \bar k$, 
that is impossible.
Therefore, $\bar k-\bar \lambda-1<\bar k/2+1$, i.e. $\bar k<2(\bar \lambda+2)$ holds. 
Substituting the expressions for $\bar k$ and $\bar \lambda$ into the previous inequality 
yields $rs<2(s+1)$. By Lemma \ref{L21}, we have $s+1\le r$.
Hence, $s+1\le r<2(s+1)/s$ and this implies that $s=1$, $r\in \{2,3\}$ and
$\Delta$ is the pentagon or the Petersen graph. In both cases 
Theorem \ref{T41} follows from Lemma \ref{L32}.
\medskip

Now we may assume $\bar \mu>1$. Since $\bar \mu<\bar k$, the left side of (\ref{E44}) is at 
most $\bar k^2$. On the other hand, if $\bar k-\bar \lambda-1>\bar k/2$ holds, then
the right side of (\ref{E44}) is more than $2\bar k \bar k/2=\bar k^2$, that is impossible. 
Hence, we have $\bar k-\bar \lambda-1\le\bar k/2$, i.e. $\bar k\le 2(\bar \lambda+1)$.

Since $\bar \mu>1$, there is an integer $\alpha_1\ge 1$ such that, for a vertex $w\in \Delta$, 
the subgraph $\Delta_1(w)$ is the $\alpha_1$-clique extension of a strongly regular Terwilliger 
graph $\Sigma$ with parameters $(v_1,k_1,\lambda_1,\mu_1)$, where 
$v_1=\displaystyle{{\frac {\bar k}{\alpha_1}}}$, $k_1=\displaystyle{{\frac {\bar \lambda-(\alpha_1-1)}{\alpha_1}}}$, 
$\mu_1=\displaystyle{{\frac {\bar \mu-1}{\alpha_1}}}$. 
Then the inequality $\bar k\le 2(\bar \lambda+1)$ is equivalent to the inequality $v_1\le 2(k_1+1)$ 
and the cardinality of the vertex set of $\Sigma$ is equal to 
$$v_1=1+k_1+k_1{\frac {(k_1-\lambda_1-1)}{\mu_1}}.$$ 

Next, $v_1\le 2(k_1+1)$ implies that 
$${\frac {k_1(k_1-\lambda_1-1)}{\mu_1}}\le k_1+1,$$
so 
$$k_1-\lambda_1-1\le \mu_1(1+1/k_1)<\mu_1+1,$$
and 
\qtnl{E45}
k_1<\lambda_1+\mu_1+2.
\eqtn

If $\mu_1=1$, then, for certain $s_1,r_1$, we have $k_1=r_1s_1$ and $\lambda_1=s_1-1$.
It follows from (\ref{E45}) that $r_1s_1<s_1-1+1+2=s_1+2$, $r_1<1+2/s_1$ and $s_1=1$, $r_1=2$. 
Hence, the graph $\Delta_1(w)$ is the $\alpha_1$-clique extension of the pentagon. 
By Lemma \ref{L32}, the graph $\Delta$ is the icosahedron and diameter of $\Gamma_1(u)$ is 3, 
that is impossible because $\Gamma$ is a Terwilliger graph.  

Hence, $\mu_1>1$. Let us consider a sequence of strongly regular graphs 
$\Sigma_1=\Sigma$, $\Sigma_2,\ldots,\Sigma_h$, $h\ge 2$ such that, for an integer $\alpha_{i+1}\ge 1$, 
the local graph of a vertex in $\Sigma_i$ is the $\alpha_{i+1}$-clique extension of 
a strongly regular Terwilliger graph $\Sigma_{i+1}$ with parameters $(v_{i+1},k_{i+1},
\lambda_{i+1},\mu_{i+1})$, $1\le i<h$ and $\mu(\Sigma_h)=1$, i.e. $\Sigma_h\in {\cal F}(s_h,r_h)$
for certain $s_h,r_h$. The sequence exists by Lemma \ref{L31}.

Assuming that $s_h>1$, we may note $k_h-\lambda_h-\mu_h=r_hs_h-(s_h-1)-1=s_h(r_h-1)>1$. 
According to Lemma \ref{L33}, we have  $k_i-\lambda_i-\mu_i>1$ for all $1\le i\le h-1$, 
which is the contradiction with (\ref{E45}). Hence, $s_h=1$ and $\Sigma_h$ is a Moore strongly 
regular graph. By Lemma \ref{L32}, diameter of $\Sigma_{h-1}$ is at least 3, 
which is the contradiction that completes the proof.\wbull


\bigskip

\begin{thebibliography}{99}

\bibitem{KP}
Jack H. Koolen, Jongyook Park: Shilla distance-regular graphs // arXiv:0902.3860 [math.CO]

\bibitem{CDG}
C. D. Godsil: Geometric distance-regular covers. New Zealand J. Math. 22 (1993), 31Ц38.

\bibitem{BCN}
A.E. Brouwer, A.M. Cohen, A. Neumaier: Distance-Regular Graphs, Springer-Verlag, Berlin Heidelberg New York, 1989.

\bibitem{GM}
A.L. Gavrilyuk, A.A. Makhnev: Locally Hoffman -- Singleton Distance-Regular Graphs // {\it to appear}.

\end{thebibliography}
\end{document}